%% file: L1TV-vixie-2007.tex
\documentclass[11pt]{article}

\pdfoutput=1
\usepackage{amsmath}
\usepackage{amssymb}
\usepackage{amsopn}
\usepackage{amsthm}
\usepackage{latexsym}
\usepackage{amsfonts}
\usepackage{subfigure}
\usepackage{url}
\usepackage{graphicx}
\usepackage{color}
\newcommand{\del}{\delta}

\newcommand{\R}{\mathbf{R}}
\newcommand{\imt}[1]{{#1}^i_*}
\newcommand{\omt}[1]{{#1}^o_*}
\newcommand{\Pm}[1]{H^1(\partial_* (#1))}
\newcommand{\Pmr}[2]{H^1(\partial_* #1 \cap #2)}

\newcommand{\Per}{\mbox{Per}}

\newcommand{\argmin}{\operatorname{argmin}}

\newtheorem{thm}{Theorem}

\newtheorem{cor}{Corollary}
\newtheorem{lemma}{Lemma}

\newtheorem{condition}{Condition}

\newtheorem{defin}{Definition}
\newtheorem{lem}{Lemma}
\newtheorem{rem}{Remark}

\usepackage{fancyhdr}
\begin{document}

\title{{\bf Some properties of minimizers for the Chan-Esedo\={g}lu $L^1$TV functional }} 
\author{Kevin R. Vixie} 
\date{\today}

\maketitle

\begin{abstract}
  We present two results characterizing minimizers of the Chan-Esedo\={g}lu $L^1$TV
  functional $F(u) \equiv \int |\nabla u | dx + \lambda \int |u - f|
  dx $; $u,f:\Bbb{R}^n \rightarrow \Bbb{R}$.  If we restrict to $u =
  \chi_{\Sigma}$ and $f = \chi_{\Omega}$, $\Sigma, \Omega \in
  \Bbb{R}^n$, the $L^1$TV functional reduces to $E(\Sigma) =
  \Per(\Sigma) + \lambda |\Sigma\vartriangle \Omega |$.  We show that
  there is a minimizer $\Sigma$ such that its boundary
  $\partial\Sigma$ lies between the union of all balls of radius
  $\frac{n}{\lambda}$ contained in $\Omega$ and the corresponding
  union of $\frac{n}{\lambda}$-balls in $\Omega^c$. We also show that
  if a ball of radius $\frac{n}{\lambda} + \epsilon$ is almost
  contained in $\Omega$, a slightly smaller concentric ball can be
  added to $\Sigma$ to get another minimizer.  Finally, we comment on
  recent results Allard has obtained on $L^1$TV minimizers and how
  these relate to our results.
\end{abstract}

\section{Introduction}

Much of the work in image analysis reduces to extracting and
processing information from images. Much of that information is, in
turn, carried by shapes present in the images.  The methods for
extracting information from images range broadly over stochastic,
wavelet, and variational or PDE based methods. In the past five to ten
years, the variational and related PDE methods have drawn a great deal
of attention.

In this paper we study one of these variational methods from a shape
processing perspective. More specifically, we establish new results
concerning the properties of exact minimizers for the rather new
Chan-Esedo\={g}lu $L^1$TV functional.  While this functional (which we
now abreviate as simply the $L^1$TV functional) is applicable to
scalar functions on $\Bbb{R}^n$, we study the functional specialized
to binary functions, i.e. binary images or shapes.

The minimization of the $L^1$TV functional,
\begin{equation}
  \label{eq:L1TV}
  u^* =  \argmin \int |\nabla u| dx + \lambda \int |u-f| dx,
\end{equation}
yields \emph{denoised images} $u$ that are smoothed yet close, in an
$L^1$ sense, to the \emph{measures image} $f$ (sometimes called the
\emph{input image} or \emph{noisy measurement}).  As is well known
from studies of the Rudin-Osher-Fatemi total
variation functional~\cite{rudin-1992-1},
\begin{equation}
  \label{eq:ROF}
  u^* =  \argmin \int |\nabla u| dx + \lambda \int |u-f|^2 dx,
\end{equation}
the total variation term $\int |\nabla u|
dx$ reduces oscillations while permitting sharp edges, something that
previous methods could not do or did very poorly. 

The change of the \emph{data fidelity term} ($\int |u-f|^2 dx$
in~(\ref{eq:ROF})) to the $L^1$ term in the $L^1$TV functional has the
effect of making that functional much more natural from a geometric
point of view.

The $L^1$TV functional was studied very carefully in a
paper by Chan and Esedoglu~\cite{chan-2005-1}. (The discrete analog of
the $L^1$TV functional had been previously studied by
Alliney~\cite{alliney-1997-1} and Nikolova~\cite{nikolova-2003-1}.)
Chan and Esedoglu~\cite{chan-2005-1} show that for binary input
images, there are also binary minimizers.  More precisely, given any
minimizer to~(\ref{eq:L1TV}) with binary input, almost every
super-levelset is the support of a binary minimizer of the same
functional. For binary input $\chi_{\Omega}$, the functional can
therefore be written as
\begin{equation}
  \label{eq:L1TVset}
  \Sigma^* = \argmin E(\Sigma) \equiv \Per(\Sigma) + \lambda|\Sigma\vartriangle\Omega|,   
\end{equation}
where $\Sigma$, $\Sigma^*$ and $\Omega$ are the supports of the binary
functions under study. Allard has recently submitted a
paper~\cite{allard-2007-1} in which he uses very intricate geometric
measure theory techniques to prove precise regularity results for
minimizers of a class of functionals which includes the $L^1$TV
functional. We comment a bit more on Allard's work in the final
section of the paper.

Our results for minimizers of the $L^1$TV functional can be viewed
results on the regularization of noisy shapes. The first result gives
us a characterization of minimizers for the case in which the noise
expresses itself as perturbations of the boundary. The second result
characterizes the $L^1$TV regularization of a binary images with
measurement noise. In discrete images this corresponds to pixels
flipping from 0 to 1 or 1 to 0 as driven by the noise process.

Now a brief outline of the paper. In the next section we present the
results for the case of $\Sigma$, $\Sigma^*$ and $\Omega$ in
$\Bbb{R}^2$. This is the case most relevant for typical images. In the
Section~\ref{sec:mtb} we prepare for the proof of these results by
introducing, in some detail, the notion of measure theoretic boundary,
exterior and interior. This permits us to avoid the intricacies of the
notion of reduced boundary. Next we prove the results for sets in
$\Bbb{R}^2$ (Section~\ref{sec:proof}).  This section is the longest
and most involved. In Section~\ref{sec:casen}, we state the theorems
for the case $n>2$ noting a few modifications that must be made.
Since all the hard parts of the proof for $n > 2$ are contained in the
$n=2$ case, we do not present the proof details.  We close
(Section~\ref{sec:disc}) with a brief discussion of our results and
their relation to one of Allard's results.

In what follows we represent minimizers of the $L^1$TV
functional~(\ref{eq:L1TVset}) by $\Sigma$, dropping the superscript
$*$ used above.

\section{Main Results (n = 2)}
\label{sec:results}

The two main results of this paper can be stated informally as
follows.  Define $R\equiv 2/\lambda$. For any $\epsilon_1, \epsilon_2
> 0$,
\begin{enumerate}
\item[(1)] any ball of radius $R$ completely contained
  in $\Omega$ is also contained in $\Sigma$, and
\item[(2)] if a ball of radius $R + \epsilon_1$ is almost contained in
  $\Omega$, then a concentric ball of radius $R - \epsilon_2$ is
  completely contained in $\Sigma$.
\end{enumerate}
More precisely we have,
\begin{thm}
\label{th:include}
Let $\Omega$ be a bounded, measurable subset of $\R^2$.  Let $\Sigma$
be any solution of (\ref{eq:L1TVset}).  Assume that a ball $B_R$ of
radius $R$ lies completely in $\Omega$: $B_R \subset \Omega$ . Then
$B_R \cup \Sigma$ is also a minimizer.  Moreover, if $B_R\subset\Omega^c$,
then $\left(B_R\cup\Sigma^c\right)^c$ is also a minimizer.
\end{thm}
and, 
\begin{thm}
\label{th:shrink}
Given $\hat{r} \in (R,\frac{\sqrt{7}}{2} R )$ and $\epsilon \in (0,1 -
\frac{1}{\sqrt{2}} )$, we can choose $\delta =
\delta(R,\hat{r},\epsilon) > 0$ such that
  \begin{equation}
    |B_{\hat{r}}\setminus\Omega| < \delta \Rightarrow
    B_{(1-\epsilon)R} \subset \Sigma .
  \end{equation}
\end{thm}
\begin{rem}
  Obvious analogs of these theorems hold in $\Bbb{R}^n$ with
  modifications commented on in Section~\ref{sec:casen}.
\end{rem}
\begin{rem}
  Theorem \ref{th:include} and the lower semicontinuity of the L1TV
  functional implies that there is a minimizer $\Sigma$ such that
  $\bigcup \{B_{\frac{2}{\lambda}}(x)\subset\Omega\} \subset \Sigma$
  and $\bigcup \{B_{\frac{2}{\lambda}}(x)\subset\Omega^c\} \subset
  \Sigma^c$. 
\end{rem}

\begin{rem}
  These theorems are close to optimal since the minimizer for $\Omega =
  B_{\frac{2}{\lambda} - \eta}$ for arbitrarily small $\eta > 0$ has unique
  minimizer $\Sigma = \emptyset$.
\end{rem}

\section{Measure Theoretic Boundary}
\label{sec:mtb}
To simplify our analysis of the energy $E(\Sigma)\equiv \Per(\Sigma) +
\lambda|\Omega\vartriangle\Sigma|$, we introduce measure theoretic
boundary, interior, and exterior.  

Define $\Per(\Sigma) \equiv \int |\nabla \chi_{\Sigma}| dx$. We say a
set in $\Bbb{R}^n$ is a \emph{set of finite perimeter} if
$\Per(\Sigma) < \infty$.  The structure theorem for sets of finite
perimeter tells us that $\Per(\Sigma) = H^{n-1}(\partial^* \Sigma)$,
where $\partial^* \Sigma$ is the \emph{reduced boundary} of $\Sigma$.
The reduced boundary is rather complicated to define and difficult to
manipulate. Instead, we use another theorem which asserts $\partial^*
\Sigma \subset \partial_* \Sigma$ and $H^{n-1}(\partial_* \Sigma -
\partial^* \Sigma) = 0$ to conclude that $\Per(\Sigma) =
H^{n-1}(\partial_* \Sigma)$, where $\partial_* \Sigma$ denotes the
\emph{measure theoretic boundary} of $\Sigma$.
(See~\cite{evans-1992-1} Theorem 2, Section 5.7 and Lemma 1, Section
5.8 for more details.) We now define measure theoretic boundary,
interior, and exterior.

\begin{defin}
\label{def:mt-stuff}
A point $x\in\Bbb{R}^n$ is in $\partial_*
A$, the \emph{measure theoretic boundary} of $A$ if
\begin{equation}
  \limsup_{r\rightarrow 0} \frac{\mathcal{L}(B(x,r)\cap A)}{r^n} > 0
\end{equation}
and
\begin{equation}
  \limsup_{r\rightarrow 0} \frac{\mathcal{L}(B(x,r)\cap A^c)}{r^n} > 0.
\end{equation}
A point $x\in\Bbb{R}^n$ is in
$A_*^i$, the \emph{measure theoretic interior} of $A$ if
\begin{equation}
  \limsup_{r\rightarrow 0} \frac{\mathcal{L}(B(x,r)\cap A^c)}{r^n} = 0.
\end{equation}
while $x\in\Bbb{R}^n$ is in
$A_*^o$, the \emph{measure theoretic exterior} of $A$ if
\begin{equation}
  \limsup_{r\rightarrow 0} \frac{\mathcal{L}(B(x,r)\cap A)}{r^n} = 0
\end{equation}
\end{defin}

\begin{lem} Let $A$ be a subset of $\Bbb{R}^n$ with finite perimeter. Then
  \begin{enumerate}
  \item $\Per(A) = H^{n-1}(\partial_* A)$
  \item $R^n = (A_*^o)\cup(\partial_* A)\cup(A_*^i)$ and the three sets
  are pairwise disjoint.
  \end{enumerate}
\end{lem}

{\noindent\bf Proof:} (1) As stated above this follows from
~\cite{evans-1992-1} Theorem 2, Section 5.7 and Lemma 1, Section 5.8.
(2) This follows directly from Definition~\ref{def:mt-stuff}.
$\blacksquare$ \bigskip

\begin{lem}
  \label{lem:mtf} Suppose $A$ and $B$ be subsets of $\Bbb{R}^n$. Then:
  \begin{enumerate}
  \item if $x\in A_*^i$ or $x\in B_*^i$ then $x\in (A\cup B)_*^i$.
  \item $\partial_* (A\cup B) \subset \partial_* A \cup \partial_* B$
  \item $\partial_* A \cup \partial_* B = (\partial_* A \cap
    B_*^i)\cup(\partial_* A \cap B_*^o)\cup(\partial_* B \cap
    A_*^i)\cup(\partial_* B \cap A_*^o)\cup(\partial_* A \cap
    \partial_* B) $
  \item (1-3) immediately imply that $ \partial_* (A\cup B) \subset
    (\partial_* A \cap B_*^o)\cup(\partial_* B \cap
    A_*^o)\cup(\partial_* A \cap \partial_* B) $
  \item $ (\partial_* A \cap B_*^o)\cup(\partial_* B \cap A_*^o)
    \subset \partial_* (A\cup B)$
  \item $\partial_* A^c = \partial_* A$.
  \item $ (A^c)_*^o = A_*^i$.
  \end{enumerate}
\end{lem}

{\noindent\bf Proof:} The lemma follows in a straightforward manner
from the definitions of measure theoretic boundary, interior and
exterior.$\blacksquare$ \bigskip

\begin{cor}
  If $H^{n-1}(\partial_* A \cap \partial_* B) = 0$ then
  \begin{enumerate}
  \item $H^{n-1}(\partial_* (A\cup B)) = H^{n-1}(\partial_* A \cap B_*^o) +
  H^{n-1}(\partial_* B \cap A_*^o)$
  \item $H^{n-1}(\partial_* (A\cap B)) = H^{n-1}(\partial_* A \cap B_*^i) +
  H^{n-1}(\partial_* B \cap A_*^i)$
  \end{enumerate}
\end{cor}

{\noindent\bf Proof:}  (1): Lemma \ref{lem:mtf}:(4-5) imply that 
\begin{eqnarray}
  H^{n-1}(\partial_* A \cap B_*^o) &+& H^{n-1}(\partial_* B \cap A_*^o) \\
  &\leq & H^{n-1}(\partial_*(A\cup B)) \\
  &\leq &  H^{n-1}(\partial_* A \cap B_*^o) + H^{n-1}(\partial_*
  B \cap A_*^o) + H^{n-1}(\partial_* A \cap \partial_* B)  
\end{eqnarray}
and the conclusion follows.  (2): This follows from (1), Lemma
\ref{lem:mtf}:(6)-(7) and the fact that $A\cap B = (A^c\cup B^c)^c$ .

$\blacksquare$ \bigskip

\begin{rem}
  Since $\partial_* A = (\partial_* A \cap B_*^i) \cup (\partial_* A \cap
  \partial_* B) \cup (\partial_* A \cap B_*^o)$, the assumption that
  $H^{n-1}(\partial_* A \cap \partial_* B) = 0$ means we can, for the
  sake of measurement, consider $\partial_* A = (\partial_* A \cap
  B_*^i) \cup (\partial_* A \cap B_*^o)$.
\end{rem}

\begin{rem}
  Now suppose that $H^{n-1}(\partial_* A) < \infty$ and $B_r$ is the
ball of radius $r$ centered at $x \in \Bbb{R}^n$ (we suppress the $x$). Note that
$\partial_* B_r = \partial B_r$. By the
coarea formula, the set of $r$'s such that $H^{n-1}(\partial_*
B_r \cap \partial_* A) > 0$ is at most countable. We conclude that
the $r$'s for which $H^{n-1}(\partial_* B_r \cap \partial_* A) = 0$
are dense and have full measure in $\Bbb{R}$. \emph{For the rest of
  this section we assume that we have chosen $r$ such that
  $H^{n-1}(\partial_* A \cap \partial_* B_r) = 0$.}
\end{rem}

\begin{thm} 
  Suppose $B_r \subset \Omega$. Define $E(\Sigma) \equiv \int |\nabla
  \chi_{\Sigma}| dx + \lambda \int |\xi_{\Sigma} -\xi_{\Omega}| dx =
  \Per(\Sigma) + \lambda |\Sigma \triangle \Omega|$.  Then
\[\Delta E = E(\Sigma\cup B_r) - E(\Sigma) = -H^{n-1}(\partial_* \Sigma \cap
  (B_r)_*^i) + H^{n-1}(\partial_* B_r \cap \Sigma_*^o) - \lambda
  |B_r\setminus \Sigma|.\]
\end{thm}

{\noindent\bf Proof:} Since $\Per(\Sigma\cup B_r) =
H^{n-1}(\partial_* \Sigma \cap (B_r)_*^o) + H^{n-1}(\partial_* B_r
\cap \Sigma_*^o)$ and $\Per(\Sigma) = H^{n-1}(\partial_* \Sigma \cap
(B_r)_*^i) + H^{n-1}(\partial_* \Sigma \cap (B_r)_*^o)$ we get
$\Per(\Sigma\cup (B_r)) - \Per(\Sigma)= -H^{n-1}(\partial_* \Sigma
\cap (B_r)_*^i) + H^{n-1}(\partial_* B_r \cap \Sigma_*^o)$.
Noting that $B_r \subset \Omega$ implies $|(\Sigma\cup B_r)
\triangle \Omega| - |\Sigma \triangle \Omega| = |B_r\setminus
\Sigma|$ finishes the proof.  $\blacksquare$ \bigskip

\begin{figure}[htbp!]
\centering
\includegraphics[width=0.65\textwidth]{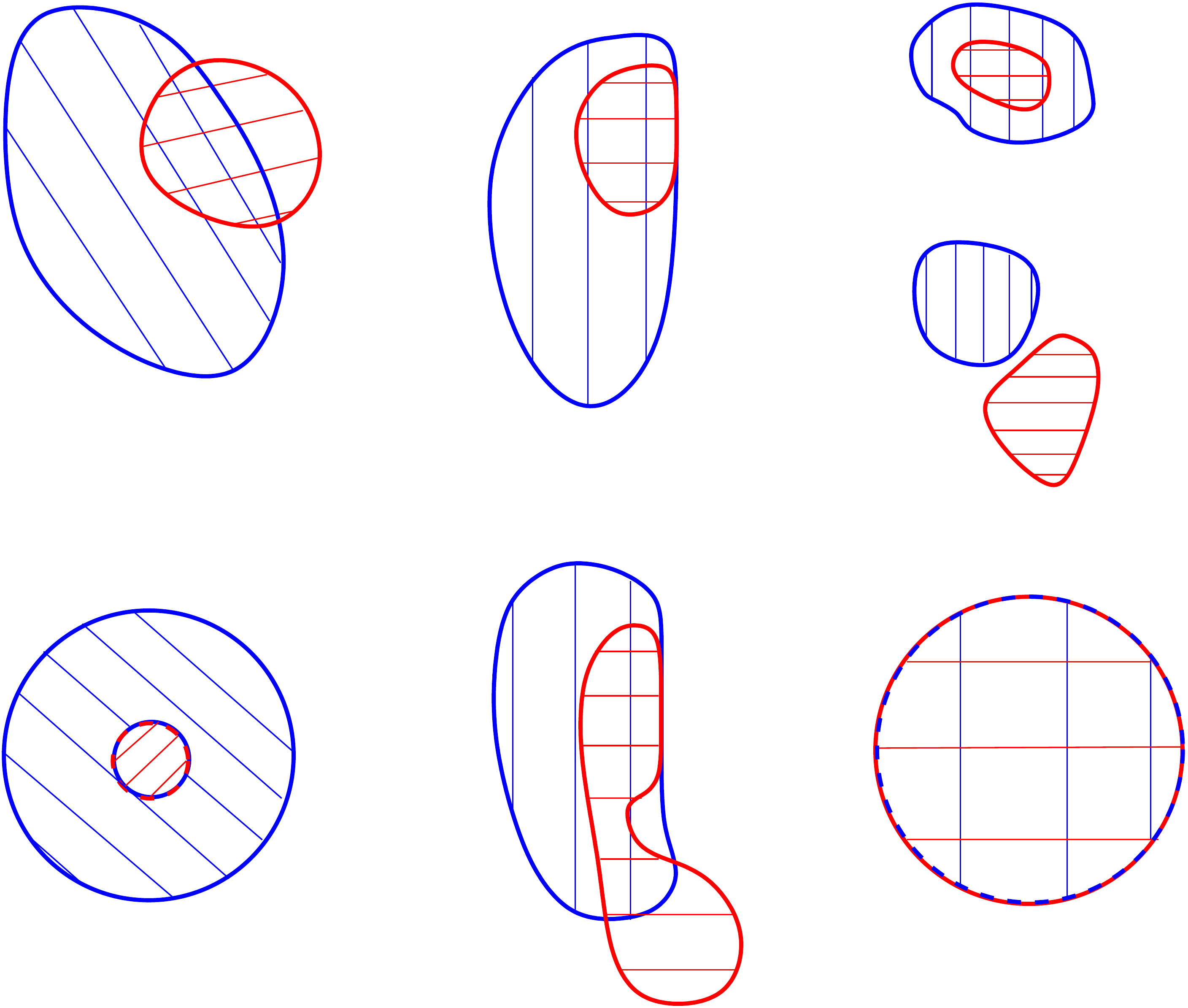} 
\caption{\label{fig:mt-illustrate} Illustration of the cases one needs
  to consider in order to understand how $\partial_* (A\cup B)$
  relates to $\partial_* A$ and $\partial_* B$}
\end{figure}

The example $A$'s and $B$'s in
Figure~\ref{fig:mt-illustrate} illustrate why the above care
is necessary.

\section{The Comparisons: proofs of Theorems~\ref{th:include} and~\ref{th:shrink}}
\label{sec:proof}

Recall that $\Sigma$ denotes a minimizer of~\ref{eq:L1TVset}.  If
$\Delta E \equiv E(B_r\cup\Sigma) - E(\Sigma) \leq 0$ then
$B_r\cup\Sigma$ must also be a minimizer.

\bigskip

\noindent {\bf Proof of Theorem~\ref{th:include}:}
Computing $\Delta E$ for $B_r\subset\Omega$ we get (for
all but countably many r):
\begin{eqnarray}
  \Delta E =& -\Pmr{\Sigma}{\imt{B_r}} + \Pmr{B_r}{\omt{\Sigma}} - 
               \lambda |B_r\setminus\Sigma| \\
           =&  - \Pmr{\Sigma}{\imt{B_r}} -\Pmr{B_r}{\imt{\Sigma}} 
               + \Pmr{B_r}{\omt{\Sigma}} + \Pmr{B_r}{\imt{\Sigma}}\nonumber \\ 
    &\;\;\;\;\;\; - \lambda |B_r\setminus\Sigma|
       - \lambda |B_r\cap\Sigma| 
               + \lambda|B_r\cap\Sigma|\\
           =&  - \Pm{B_r\cap\Sigma} + \Pm{B_r}    
               -\lambda |B_r| + \lambda |B_r\cap\Sigma|\\
           =&  \left(\Pm{B_r} - \lambda |B_r|\right) 
             + \left(\lambda |B_r\cap\Sigma| - \Pm{B_r\cap\Sigma}\right).\label{eq:thm1prf}
\end{eqnarray}
This is illustrated in Figure~\ref{fig:B-in-Omega} below.
\begin{figure}[htbp!]
  
\begin{center}
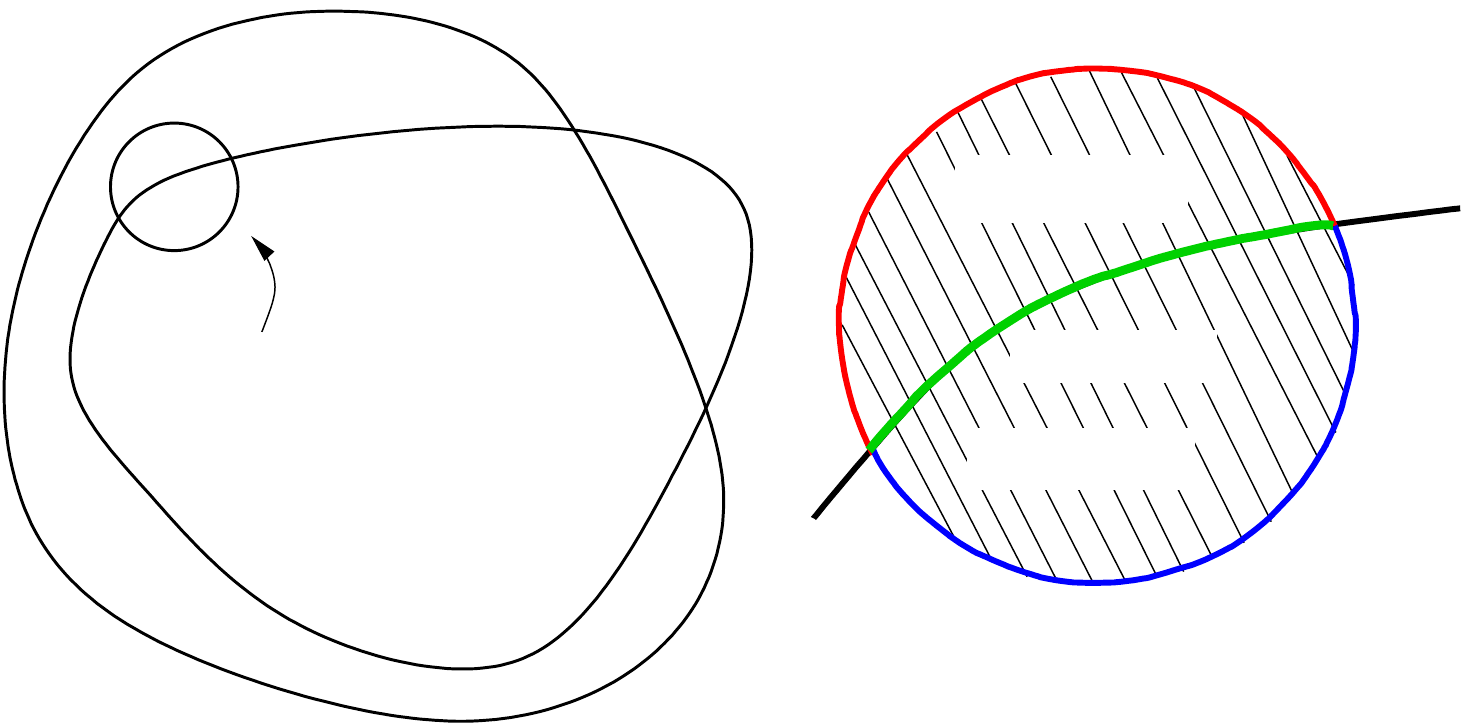
\caption{\label{fig:B-in-Omega}An illustration useful for computing $\Delta E$ when
  $B_r\cap\Omega^c = \emptyset$.}
\end{center}  
\end{figure}

Now we choose a sequence of radii $r_i < R$ converging to $R$, for
which~(\ref{eq:thm1prf}) holds.

Defining $\rho_i$ by $\pi\rho_i^2 = |B_{r_i}\cap\Sigma|$, $\rho_i^*$
by $2\pi\rho_i^* = \Pm{B_{r_i}\cap\Sigma}$ and remembering that $R
\equiv \frac{2}{\lambda}$, we have that
\begin{eqnarray}
  \Delta E &= \left(2\pi r_i - \frac{2}{R}\pi r_i^2 \right) 
             + \left( \frac{2}{R}\pi\rho_i^2 - 2\pi\rho_i^* \right)\\
           &= 2\pi r_i (1-\frac{r_i}{R}) + 2\pi\rho_i (\frac{\rho_i}{R} 
             - \frac{\rho_i^*}{\rho_i}). \label{eq:thm1ineq}
\end{eqnarray}
Note that the isoperimetric inequality gives $\frac{\rho_i^*}{\rho_i}
\geq 1$ for all $i$, that $\frac{\rho_i}{R} < 1$ for all $i$ and that
$(1-\frac{r_i}{R})\rightarrow_{i\rightarrow\infty} 0$.  The right hand
side of~(\ref{eq:thm1ineq}) converges therefore to zero.  Using the
fact that $\Delta E$ is lower semicontinuous for sequences in $L^1$,
(which follows from the lower semicontinuity of the BV seminorm), we
conclude that $\Delta E (B_R) \leq 0$. We conclude that
$\Sigma\cup B_R$ is also a minimizer.

Finally, we note that $E_{\Omega}(\Sigma) \equiv \Per(\Sigma) +
\lambda|\Sigma\vartriangle\Omega| = E_{\Omega^c}(\Sigma^c) \equiv
\Per(\Sigma^c) + \lambda|\Sigma^c\vartriangle\Omega^c|$.  From this we
deduce that $\Sigma$ minimizes $E_{\Omega} \Leftrightarrow \Sigma^c$
minimizes $E_{\Omega^c}$. Therefore, $B_R\subset\Omega^c$ implies
$(\Sigma^c \cup B_R)^c$ is also a minimizer.  $\blacksquare$

\bigskip

\noindent {\bf Proof of Theorem~\ref{th:shrink}:}
In the case that $B_r\cap\Omega^c \neq \emptyset$,
\begin{eqnarray}
  \Delta E &=& -\Pmr{\Sigma}{\imt{B_r}} + \Pmr{B_r}{\omt{\Sigma}})  
               - \lambda |B_r\setminus\Sigma| + 
                 2\lambda |B_r\cap\Omega^c\cap\Sigma^c |\\
           &\leq& -\Pmr{\Sigma}{\imt{B_r}} + \Pmr{B_r}{\omt{\Sigma}}  
                 - \lambda |B_r\setminus\Sigma| + 
                 2\lambda |B_r\cap\Omega^c|\\
            &=& -\Pmr{\Sigma}{\imt{B_r}} + \Pmr{B_r}{\omt{\Sigma}}  
                 - \lambda |B_r\setminus\Sigma| + 
                 2\lambda |B_r\setminus\Omega|\\
           &=&  \left(\Pm{B_r} - \lambda |B_r|\right) 
             + \left(\lambda |B_r\cap\Sigma| - \Pm{B_r\cap\Sigma}\right)
             + 2\lambda |B_r\setminus\Omega|\\
           &=& \left(\Per(B_r) - \lambda |B_r|\right) 
             + \left(\lambda |B_r\cap\Sigma| - \Per(B_r\cap\Sigma)\right)
             + 2\lambda |B_r\setminus\Omega|. \label{eq:lemma}
\end{eqnarray}

This is illustrated in Figure~\ref{fig:B-in-Omega-delta} below.
\begin{figure}[htbp!]
\begin{center}
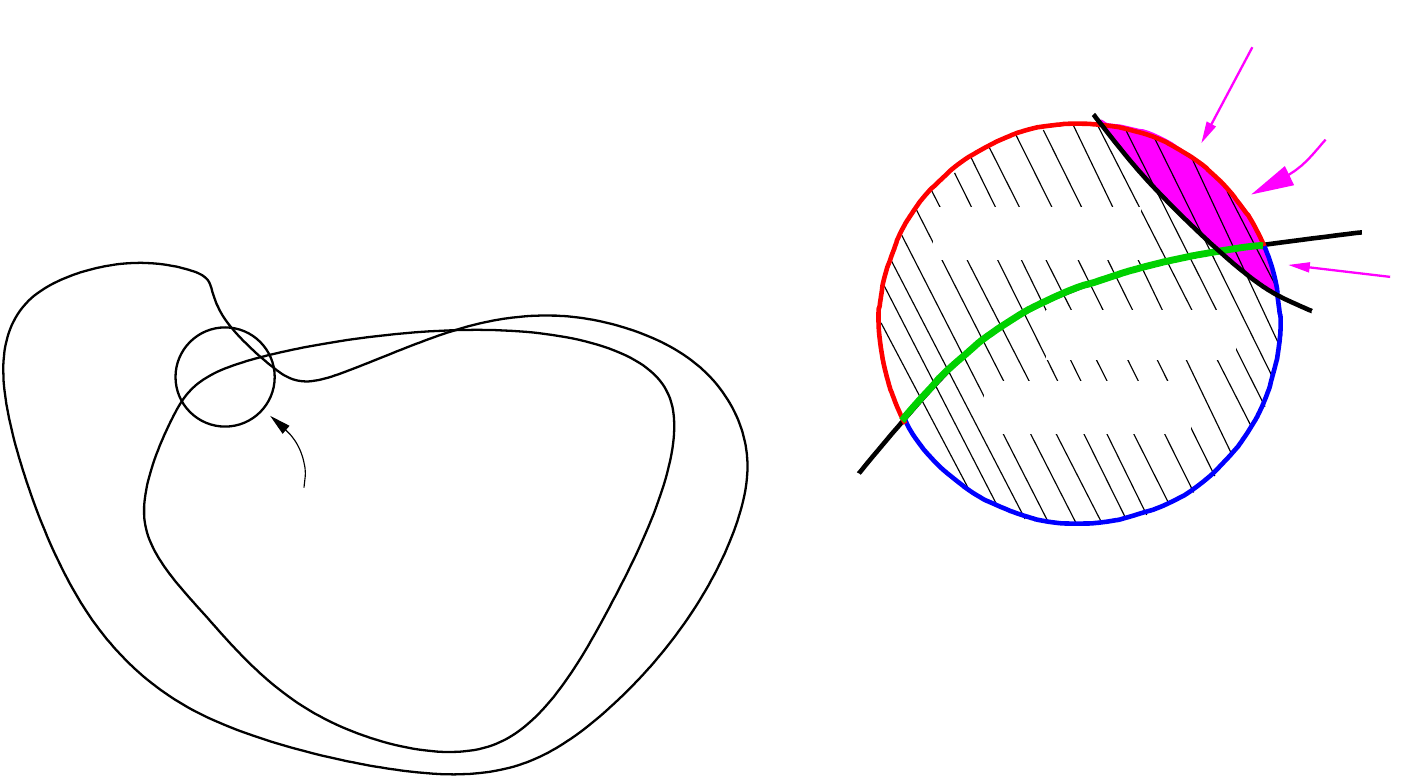
\caption{ \label{fig:B-in-Omega-delta}An illustration useful for
  computing $\Delta E$ when $B_r\cap\Omega^c \neq \emptyset$}
\end{center}  
\end{figure}

\begin{lem}
  $|B_R\setminus\Sigma| \leq 6 \delta$   
\end{lem}
\noindent {\bf Proof:} Since we assume $\Sigma$ is a minimizer,
$\Delta E \geq 0$.  We will perturb with balls of radius $r\leq
\hat{r}$.  Then, $|B_r\setminus\Omega| \leq
|B_{\hat{r}}\setminus\Omega| := \del$. These assumptions together
with~(\ref{eq:lemma}) and the isoperimetric inequality
($\Per(B_r\cap\Sigma) \geq 2\sqrt{\pi}|B_r\cap\Sigma|^{\frac{1}{2}}$)
imply:
\begin{eqnarray}
  0 & \leq & \Delta E \leq 2\pi r - \lambda \pi r^2 + \lambda |B_r\cap\Sigma| - 2\sqrt{\pi}|B_r\cap\Sigma|^{\frac{1}{2}} + 2\lambda\del. \\
    &  =  &  \lambda|B_r\cap\Sigma| - 2\sqrt{\pi}|B_r\cap\Sigma|^{\frac{1}{2}} + \Big( 2\lambda\del + 2\pi r - \lambda\pi r^2 \Big)\\
    &  =  &  f(\xi) \equiv \lambda\xi^2 - 2\sqrt{\pi}\xi + \Big( 2\lambda\del + 2\pi r - \lambda\pi r^2 \Big), \text{ ($\xi\equiv |B_r\cap\Sigma|^{\frac{1}{2}}$) }\\
    &  =  & \frac{2}{R}\xi^2 - 2\sqrt{\pi}\xi + \Big( \frac{4\del}{R} + 2\pi r - \frac{2\pi r^2}{R} \Big) \text{ (recalling $R=\frac{2}{\lambda}$) \label{eq:poly2}}.
\end{eqnarray}
In view of the last inequality, we describe values of $\xi$ for which
$f(\xi)\geq 0$.  For a given $r$, the zeros of $f(\xi)$ are at:
\begin{equation}
\label{eq:genr}
\xi_\pm(r) = \sqrt{\frac{\pi R^2}{4}} \; \pm \; \sqrt{ \frac{\pi R^2}{4} + (\pi r (r-R)) - 2\del}.
\end{equation}
Thus, for all $r\leq\hat{r}$, we have:
\begin{equation}
  \mbox{Either } |B_r\cap\Sigma|^{\frac{1}{2}} \leq \xi_-(r) \mbox{, or } |B_r\cap\Sigma|^{\frac{1}{2}} \geq \xi_+(r).
\end{equation}
If we take $r = \hat{r} > R $, then $2 \pi r(r - R) > 0$, and assuming
\begin{condition}
\label{cn:first}
  \begin{equation}
    \del < \frac{\pi \hat{r}}{2} (\hat{r}-R)   
  \end{equation}
\end{condition}
\noindent implies $\xi_{-}(r) < 0 $. This implies
\begin{equation}
  |B_{\hat{r}} \cap \Sigma | \geq \xi_+^2(\hat{r}) > \pi R^2.
\end{equation}
Since $ r = \hat{r} < \frac{\sqrt{7}}{2} R$ we get
\begin{equation}
\label{eq:br1}
|B_R \cap \Sigma | > |B_r \cap \Sigma | -
\frac{3\pi R^2}{4} > \frac{\pi R^2}{4}  
\end{equation}
Now we consider $\xi_\pm(R)$:
\begin{equation}
\label{eq:br2} 
  \xi_\pm(R) = \sqrt{\frac{\pi R^2}{4}} \pm \sqrt{\frac{\pi R^2}{4} - 2 \del}.
\end{equation}
Assuming
\begin{condition}
  \label{cn:cond-1-8}
  \begin{equation}
        \del < \frac{\pi R^2}{8},    
  \end{equation}
\end{condition}
we get that $\xi_\pm(R)$ are real and distinct. Since
\begin{equation}
  \label{eq:br3}
\xi_{-}^2(R) < \frac{\pi R^2}{4} < |B_R \cap\Sigma|,  
\end{equation}
we conclude that
 \begin{equation}
\label{eq:br4}
   \xi^2_+ \leq |B_R \cap \Sigma |.
 \end{equation}
Computing, we get
\begin{eqnarray}
  \label{eq:xibegin}
  \xi^2_+ &=& \left(\sqrt{\frac{\pi R^2}{4}} + 
                  \sqrt{\frac{\pi R^2}{4} - 2 \del}\right)^2 \\
          &=& \frac{\pi R^2}{2} - 2\delta + \frac{\pi R^2}{2}
                 \sqrt{1 - \frac{8\delta}{\pi R^2}}\\
          &\geq& \frac{\pi R^2}{2} - 2\del + \frac{\pi R^2}{2} 
                  \left(1 - \alpha\frac{8\del}{\pi R^2}\right) 
            \text{ $\left(\text{ assuming } \delta \leq \frac{\pi R^2}{8}
                               \frac{2\alpha - 1}{\alpha^2}\right)$ }\label{eq:alphaCN}\\
           &=& \pi R^2 - (2 + 4\alpha)\del. \label{eq:xiend}     
\end{eqnarray}
Choosing $\alpha = 1$, and noting that Condition~\ref{cn:cond-1-8}
then implies the assumption in~\ref{eq:alphaCN} is satisfied, we get
\begin{equation}
  |B_R \cap \Sigma | \geq \xi^2_+ \geq \pi R^2 - 6\del.
\end{equation}
This gives 
\begin{equation}
  |B_R \setminus \Sigma |  \leq 6\del
\end{equation}
as advertised. $\blacksquare$
\bigskip
\begin{rem} What if either $\hat{r}$ or $R$ are radii such
  that~(\ref{eq:lemma}) (and therefore~(\ref{eq:genr})) does not hold?
  We can simply choose another $\tilde{r} < \hat{r}$ arbitrarily close
  to $\hat{r}$, for which~(\ref{eq:lemma}) does hold. The
  $\tilde{\delta} \equiv |B_{\tilde{r}} \setminus \Omega |$ will be no
  greater than, and arbitrarily close to, $\delta$.  As we will see,
  the only conditions on $\delta$ that are not functions of $R$ and
  $\epsilon$ are those in Condition \ref{cn:first}. Therefore, if we
  replace Condition~\ref{cn:first} with
  \begin{equation}
    \label{cn:firstprime}
    \del < \frac{\pi R}{4} (\hat{r}-R)   
  \end{equation} 
  we know that the delta chosen for any $\hat{r}$ will permit us to
  arrive at the conclusions of this lemma, even in cases where we have
  to perturb $\hat{r}$. Next we choose a sequence of $r_i > R$
  converging monotonically to $R$ for which the inequality does work.
  Equation~(\ref{eq:br1}) is still valid if we replace $B_R$ with
  $B_{r_i}$.  Equation~(\ref{eq:br2}) can be slightly modified using~(\ref{eq:genr}) to
\begin{equation}
\label{eq:br2r} 
  \xi_\pm(r_i) = \sqrt{\frac{\pi R^2}{4}} \pm \sqrt{\frac{\pi R^2}{4} - 2 \del_i}.
\end{equation}
where the $\del_i < \del$ and $\del_i \rightarrow \del$ as
$i\rightarrow\infty$.  Now, simply repeating the derivation in lines
(\ref{eq:xibegin}) to (\ref{eq:xiend}), gives 
\begin{equation}
  |B_R\cap\Sigma| = \lim_{i\rightarrow\infty} |B_{r_i}\cap\Sigma| \geq \pi R^2 - 6\del =
     \lim_{i\rightarrow\infty} \pi R^2 - 6\del_i
\end{equation}
\end{rem}

Now we continue with the proof of Theorem~\ref{th:shrink}. Computing
(again and less optimally, but sufficiently for our purposes) the
change in energy when we add a ball $B_r$ to $\Sigma$ for $r\in(0,R)$,
we get
\begin{eqnarray}
\label{changeinenergy}
\Delta E = E(\Sigma\cup B_r) &-& E(\Sigma)\\
&\leq&  -\Pmr{\Sigma}{\imt{B_r}} + \Pmr{B_r}{\omt{\Sigma}} +
\lambda |B_r\setminus\Sigma| \\
&=& -\Per(\Sigma;B_r) + \Pmr{B_r}{\omt{\Sigma}} +
\lambda |B_r\setminus\Sigma|.
\end{eqnarray}

By the coarea formula and properties of the measure theoretic
exterior,
\begin{equation}
\label{polarcoords}
|B_r\setminus\Sigma| = |B_r\cap\Sigma_*^o| = \int_0^r \Pmr{B_{\xi}}{\omt{\Sigma}} \, d\xi.
\end{equation}
By the relative isoperimetric inequality applied in the ball
$B_r(x_0)$,
\begin{equation}
\Per(\Sigma;B_r) \geq C \min \left\{ |B_r\setminus\Sigma|^\frac{1}{2}, 
|\Sigma\cap B_r|^\frac{1}{2} \right\}.
\end{equation}
Assuming           
\begin{condition}
\label{cond:first}
$6\delta  < \frac{1}{4} \pi R^2$            
\end{condition}
implies $|B_R\setminus\Sigma| <
\frac{1}{4} \pi R^2$. Assuming $r > \frac{R}{\sqrt{2}}$ implies that
$|B_r\setminus\Sigma|\leq|\Sigma\cap B_r|$ and consequently
\begin{equation}
\label{relativeiso}
\Per(\Sigma;B_r) \geq C |B_r\setminus\Sigma|^\frac{1}{2}.
\end{equation}
This gives a
condition on $\epsilon$:
\begin{condition}
  $\epsilon < 1 - \frac{1}{\sqrt{2}}$.
\end{condition}

Define $v(r):=|B_r\setminus\Sigma|$.  By differentiating
(\ref{polarcoords}) with respect to $r$, and using (\ref{relativeiso})
we see that the inequality concerning the change in energy given in
(\ref{changeinenergy}) can be written as
\begin{equation}
\label{diffineq}
E(\Sigma\cup B_r) - E(\Sigma) \leq \lambda v(r) - C\sqrt{v(r)} + v'(r).
\end{equation}
We will use the differential expression on the right to show that the
change in energy on the left has to be negative for some $r$ close to
$R$.  

\begin{rem}
  Note that by choosing $\delta$ small enough, we can make $v(r)$
  arbitrarily small and obtain $\lambda v(r) - C\sqrt{v(r)} < 0$; if
  the right hand side is positive then we have $v'(r) > 0$.  This in
  turn means that $v(r)$ decreases as r gets smaller.  We exploit this
  to force the right hand side to zero.
\end{rem}
\medskip

\begin{lemma}
  $v'(r) - C\sqrt{v(r)} + \lambda v(r) \leq 0$ for a set of $r \in
  \left( (1-\epsilon)R,R\right)$ with positive measure.
\end{lemma}

\smallskip

\noindent {\em Proof of lemma:} Assume
\begin{equation}
\label{diffineq2}
v'(r) - C\sqrt{v(r)} + \lambda v(r) \geq 0 \mbox{ for a.e. } r \in \left( (1-\epsilon)R,R\right),
\end{equation}
otherwise we are done.
Let $w(s):=e^{-\lambda s}v(R-s)$.
Then (\ref{diffineq2}) turns into
\begin{equation}
\label{eq:second-diff}
  w'(s) + Ce^\frac{-\lambda s}{2} \sqrt{w(s)} \leq 0 \mbox{ for a.e. } s \in \left(0,\epsilon R\right).
\end{equation}
with the initial condition $w(0)=|B_R\setminus\Sigma|$ and $w(s) \geq 0$.
Solutions of this differential inequality can be bounded from above by
solutions of the following differential equality:
\begin{equation}
\label{diffeq}
\begin{split}
\Bar{w}' &= -Ce^{\frac{-\lambda s}{2}}\sqrt{\Bar{w}}.\\
\Bar{w}(0) &= |B_R\setminus\Sigma|\text{ and } \Bar{w} \geq 0.
\end{split}
\end{equation}
The solution is
\begin{eqnarray*}
\sqrt{\Bar{w}(s)} & = &\max\left( 0 , \frac{C}{\lambda} \Big(e^{-\frac{\lambda}{2}s}-1\Big) 
+ \sqrt{|B_R\setminus\Sigma|}\right)\\
                  & = & \max\left( 0 , \frac{CR}{2} \Big(e^{-\frac{s}{R}}-1\Big)
                   + \sqrt{|B_R\setminus\Sigma|}\right).
\end{eqnarray*}
Therefore if $|B_R\setminus\Sigma|\leq 6\delta$ and 
\begin{condition}
  $6\delta \leq \alpha$, where $\alpha$ is any solution to
  \begin{equation}
    \frac{CR}{2} \Big(e^{-\frac{\epsilon R}{R}}-1\Big) +
  \sqrt{\alpha} = \frac{CR}{2}\Big(e^{-\epsilon} - 1 \Big) +
  \sqrt{\alpha} < 0
  \end{equation}
i.e., we have 
\begin{equation}
 \delta < \frac{C^2 R^2}{24}( 1 - e^{-\epsilon})^2,
\end{equation}
\end{condition}
\noindent then we have a set of $r$ with positive measure in $((1-\epsilon)R,R)$
such that $v(r) = 0$ and $v'(r) = 0$. $\blacksquare$

This lemma immediately implies that for some $r\in ((1-\epsilon)R,R)$, 
$B_{r}\cup\Sigma$ is also a minimizer. $\blacksquare$

\bigskip

\section{\bf The Case  $n > 2$}
\label{sec:casen}
The analogs for Theorems~\ref{th:include} and~\ref{th:shrink} in $\Bbb{R}^n$ are:
\begin{thm}
\label{th:include-n}
Let $\Omega$ be a bounded, measurable subset of $\R^n$.  Let $\Sigma$
be any solution of (\ref{eq:L1TVset}).  Assume that a ball $B_R$, $R =
\frac{n}{\lambda}$ of radius $R$ lies completely in $\Omega$: $B_R
\subset \Omega$ . Then $B_R \cup \Sigma$ is also a minimizer.
Moreover, if $B_R\subset\Omega^c$, then
$\left(B_R\cup\Sigma^c\right)^c$ is also a minimizer.
\end{thm}
\begin{thm}
\label{th:shrink-n}
Given 
\begin{equation*}
\hat{r} \in \left(R,\left(2 - \left(\frac{n-1}{n}\right)^n\right)^{\frac{1}{n}} R \right)  
\end{equation*}
and 
\begin{equation*}
\epsilon \in \left(0,1 - \frac{1}{2^{\frac{1}{n}}} \right)  
\end{equation*}
we can choose $\delta > 0$ such that
  \begin{equation}
    |B_{\hat{r}}\setminus\Omega| < \delta \Rightarrow
    B_{(1-\epsilon)R} \cup \Sigma \text{ is also a minimizer. }
  \end{equation}
\end{thm}
We do not present the proofs, since they are very similar to the $n=2$
case. In particular, making the replacement $R = \frac{2}{\lambda}
\rightarrow R = \frac{n}{\lambda}$ enables us to use the proof of
Theorem~\ref{th:include}, with obvious modifications, to obtain
Theorem~\ref{th:include-n}. Likewise, we can use the proof of
Theorem~\ref{th:shrink} to prove Theorem~\ref{th:shrink-n}, with
modifications noted below.
\begin{itemize}
\item[(1)] Again, $R = \frac{2}{\lambda} \rightarrow R =
  \frac{n}{\lambda}$,
\item[(2)] We define $\xi \equiv |B_r\cap\Sigma|^{\frac{1}{n}}$. Let
  $\alpha_n$ be the volume of the ball with unit radius in
  $\Bbb{R}^n$. The polynomial in~(\ref{eq:poly2}) then gets
  replaced by
  \begin{equation}
    \label{eq:polyn}
    \frac{n}{R}\xi^n - n \alpha_n^{\frac{1}{n}}\xi^{n-1} + \left(\frac{2n\delta}{R} 
     + n \alpha_n r^{n-1} - \frac{n\alpha_n r^n}{R}  \right).
  \end{equation}
  Since we are interested in the roots of this polynomial, we look at
\begin{equation}
  \label{eq:polyn-clean}
  \xi^{n-1}\left( \alpha_n^{\frac{1}{n}} R - \xi \right) = 2\delta + 
  \alpha_n r^{n-1}\left( R - r\right).
\end{equation}  
\item[(3)] We replace the right hand side of~(\ref{diffineq}) with
  \begin{equation}
    \label{eq:diffineq-n}
    \lambda v(r) - C v(r)^{\frac{n-1}{n}} + v'(r)
  \end{equation} 
  which gives us
\begin{equation}
  \label{eq:second-diff-n}
  w'(s) + C e^{\frac{-\lambda s}{n}}\left(w(s)\right)^{\frac{n-1}{n}} \leq 0
\end{equation}
in place of~(\ref{eq:second-diff}).
\end{itemize}

\section{Discussion}
\label{sec:disc}

As mentioned in the introduction, Allard~\cite{allard-2007-1} has
recently produced an extensive study of the regularity of minimizers
for a class of functionals including the $L^1$TV functional. In
this work he uses geometric measure theory  techniques originally
developed to address minimal surface problems. As a result, his $n$
cannot exceed $7$.  In our work we have used simpler pieces of
geometric measure theory, specifically the structure theory for sets
of finite perimeter. The weaker regularity results we use -- simply
what one gets from $\Sigma$ having finite perimeter -- are not limited
to $n \leq 7$.

Allard proves that the total mean curvature of minimizers is bounded
by $\lambda$, as suggested by a naive calculation with the formal
Euler-Lagrange equation. For spheres, this corresponds to a radius of
curvature of $\frac{n-1}{\lambda}$.  In our work we find that spheres
(balls) of radii $\frac{n}{\lambda}$ play a critical role.  Such a
sphere has a total mean curvature of $\frac{n - 1}{n}\lambda$.  This
second, bigger radius characterizes the global nature of the
minimizers.  Indeed if one can contain $\Omega$ in a ball of radius
$\frac{n}{\lambda} -\epsilon$, where $0 < \epsilon$, then the unique
solution is the empty set. This follows from a monotonicity result
proved by Yin in~\cite{yin-2006-1}.  It also follows from monotonicity
results in Allard's paper~\cite{allard-2007-1}.

Another previous work that needs to be mentioned is the work of Italo
Tamanini and collaborators (see~\cite{tamanini-1996-1} and
references). Instead of using knowledge of $\Omega$ to deduce
properties of the minimizer $\Sigma$, they use weaker properties of
$\Sigma$ to establish stonger properties of the same $\Sigma$. (In
particular, if a ball of a particular radius is almost contained in
$\Sigma$ then the ball with half the radius and same center is
completely contained in $\Sigma$. This is similar to our
Thereom~\ref{th:shrink}.) These types of regularity properties of
minimizers are not very useful for computation or in the establishment
of minimizer properties based only on realistically obtainable
knowledge.

There are numerous potential directions in which to advance to these
results and the results reported in \cite{allard-2007-1,yin-2006-1}.
Generalization to anisotropic energies (see~\cite{esedoglu-2004-1} for
example), the construction of hybrid analytic-numerical algorithms for
$L^1$TV minimization, and the exploitation and analysis of the scale
decomposition properties of the $L^1$TV functional are three that come
easily to mind.

\section{Acknowledgements}
\label{sec:ackn}

It gives me pleasure to acknowledge significant conversations with
Selim Esedo\={g}lu who contributed substantial suggestions for
improvment as well as suggesting the initial idea that started this
paper. I also acknowledge the pleasant benefit of discussions with Bill
Allard. The research was made possible by funding from Los Alamos
National Laboratory and the Department of Energy. I would also like to
acknowledge the hospitality of the Institute for Pure and Applied
Analysis at UCLA where I was able to work uninterrupted on the
research reported here.

\bibliographystyle{plain}
\bibliography{/home/vixie/projects/templates/the_bib}

\end{document}

%% file: B-in-Omega-ltx.pdf_t
\begin{picture}(0,0)%
\includegraphics{B-in-Omega-ltx.pdf}%
\end{picture}%
\setlength{\unitlength}{3729sp}%
\begingroup\makeatletter\ifx\SetFigFont\undefined%
\gdef\SetFigFont#1#2#3#4#5{%
  \reset@font\fontsize{#1}{#2pt}%
  \fontfamily{#3}\fontseries{#4}\fontshape{#5}%
  \selectfont}%
\fi\endgroup%
\begin{picture}(7452,3685)(3022,-3454)
\put(8193,-1630){\makebox(0,0)[lb]{\smash{{\SetFigFont{11}{13.2}{\rmdefault}{\mddefault}{\updefault}{\color[rgb]{0,.82,0}$\partial_*\Sigma\cap{B_{r}^i}_*$}%
}}}}
\put(4552,-187){\makebox(0,0)[lb]{\smash{{\SetFigFont{11}{13.2}{\rmdefault}{\mddefault}{\updefault}{\color[rgb]{0,0,0}$\Omega$}%
}}}}
\put(6377,-1032){\makebox(0,0)[lb]{\smash{{\SetFigFont{11}{13.2}{\rmdefault}{\mddefault}{\updefault}{\color[rgb]{0,0,0}$\Sigma$}%
}}}}
\put(4117,-1672){\makebox(0,0)[lb]{\smash{{\SetFigFont{11}{13.2}{\rmdefault}{\mddefault}{\updefault}{\color[rgb]{0,0,0}$B_r$}%
}}}}
\put(3887,-907){\makebox(0,0)[lb]{\smash{{\SetFigFont{11}{13.2}{\rmdefault}{\mddefault}{\updefault}{\color[rgb]{0,0,0}$1$}%
}}}}
\put(3757,-622){\makebox(0,0)[lb]{\smash{{\SetFigFont{11}{13.2}{\rmdefault}{\mddefault}{\updefault}{\color[rgb]{0,0,0}$2$}%
}}}}
\put(8195, 60){\makebox(0,0)[lb]{\smash{{\SetFigFont{11}{13.2}{\rmdefault}{\mddefault}{\updefault}{\color[rgb]{1,0,0}$\partial_* B_r \cap \Sigma_*^o$}%
}}}}
\put(8205,-2990){\makebox(0,0)[lb]{\smash{{\SetFigFont{11}{13.2}{\rmdefault}{\mddefault}{\updefault}{\color[rgb]{0,0,1}$\partial_* B_r \cap \Sigma_*^i$}%
}}}}
\put(7917,-781){\makebox(0,0)[lb]{\smash{{\SetFigFont{11}{13.2}{\rmdefault}{\mddefault}{\updefault}{\color[rgb]{0,0,0}$2 = | B_r \cap \Sigma^c|$}%
}}}}
\put(7974,-2179){\makebox(0,0)[lb]{\smash{{\SetFigFont{11}{13.2}{\rmdefault}{\mddefault}{\updefault}{\color[rgb]{0,0,0}$1 = | B_r \cap \Sigma |$}%
}}}}
\end{picture}%

%% file: B-in-Omega-delta-ltx.pdf_t
\begin{picture}(0,0)%
\includegraphics{B-in-Omega-delta-ltx.pdf}%
\end{picture}%
\setlength{\unitlength}{2901sp}%
\begingroup\makeatletter\ifx\SetFigFont\undefined%
\gdef\SetFigFont#1#2#3#4#5{%
  \reset@font\fontsize{#1}{#2pt}%
  \fontfamily{#3}\fontseries{#4}\fontshape{#5}%
  \selectfont}%
\fi\endgroup%
\begin{picture}(9201,5058)(1545,-4389)
\put(4911,-2457){\makebox(0,0)[lb]{\smash{{\SetFigFont{8}{9.6}{\rmdefault}{\mddefault}{\updefault}{\color[rgb]{0,0,0}$\Sigma$}%
}}}}
\put(1956,-1572){\makebox(0,0)[lb]{\smash{{\SetFigFont{8}{9.6}{\rmdefault}{\mddefault}{\updefault}{\color[rgb]{0,0,0}$\Omega$}%
}}}}
\put(2886,-1677){\makebox(0,0)[lb]{\smash{{\SetFigFont{8}{9.6}{\rmdefault}{\mddefault}{\updefault}{\color[rgb]{0,0,0}$2$}%
}}}}
\put(3006,-1992){\makebox(0,0)[lb]{\smash{{\SetFigFont{8}{9.6}{\rmdefault}{\mddefault}{\updefault}{\color[rgb]{0,0,0}$1$}%
}}}}
\put(9711,498){\makebox(0,0)[lb]{\smash{{\SetFigFont{8}{9.6}{\rmdefault}{\mddefault}{\updefault}{\color[rgb]{1,0,1}$3 = | B_r \cap \Omega^c \cap \Sigma^c |$}%
}}}}
\put(10311,-237){\makebox(0,0)[lb]{\smash{{\SetFigFont{8}{9.6}{\rmdefault}{\mddefault}{\updefault}{\color[rgb]{1,0,1}$| B_r \cap \Omega^c | = 3 + 4$}%
}}}}
\put(10731,-1152){\makebox(0,0)[lb]{\smash{{\SetFigFont{8}{9.6}{\rmdefault}{\mddefault}{\updefault}{\color[rgb]{1,0,1}$4 = | B_r \cap \Omega^c \cap \Sigma |$}%
}}}}
\put(8166,-3012){\makebox(0,0)[lb]{\smash{{\SetFigFont{8}{9.6}{\rmdefault}{\mddefault}{\updefault}{\color[rgb]{0,0,1}$\partial_* B_r \cap \Sigma_*^i$}%
}}}}
\put(7866, 78){\makebox(0,0)[lb]{\smash{{\SetFigFont{8}{9.6}{\rmdefault}{\mddefault}{\updefault}{\color[rgb]{1,0,0}$\partial_* B_r \cap \Sigma_*^o$}%
}}}}
\put(8406,-1537){\makebox(0,0)[lb]{\smash{{\SetFigFont{8}{9.6}{\rmdefault}{\mddefault}{\updefault}{\color[rgb]{0,.82,0}$\partial_* \Sigma \cap B_{r*}^i$}%
}}}}
\put(7691,-872){\makebox(0,0)[lb]{\smash{{\SetFigFont{8}{9.6}{\rmdefault}{\mddefault}{\updefault}{\color[rgb]{0,0,0}$2 = | B_r \cap \Sigma^c |$}%
}}}}
\put(8031,-2017){\makebox(0,0)[lb]{\smash{{\SetFigFont{8}{9.6}{\rmdefault}{\mddefault}{\updefault}{\color[rgb]{0,0,0}$1 = | B_r \cap \Sigma |$}%
}}}}
\put(3451,-2701){\makebox(0,0)[lb]{\smash{{\SetFigFont{8}{9.6}{\rmdefault}{\mddefault}{\updefault}{\color[rgb]{0,0,0}$B_r$}%
}}}}
\end{picture}%